\documentclass[conference]{IEEEtran}
\usepackage[T1]{fontenc}
\usepackage{stackengine}
\IEEEoverridecommandlockouts
\usepackage{cite}
\usepackage{amsmath,amssymb,amsfonts, array}
\usepackage{algorithmic}
\usepackage{graphicx}
\usepackage{textcomp}
\usepackage{xcolor}
\usepackage{algorithm,algorithmic}
\usepackage{mathtools, stackengine, mathrsfs, setspace, bm, multirow, tikz, pgfplots, adjustbox, calc, url, xtab}
\usepackage[colorlinks=true, citecolor=blue, linkcolor=black!20!red]{hyperref}
\usepackage[caption=false,font=footnotesize]{subfig}
\usepackage{xparse}
\usetikzlibrary{arrows,automata}

\allowdisplaybreaks

\begin{document}
\bstctlcite{IEEEexample:BSTcontrol}
\title{Predictive Prescription of Unit Commitment Decisions Under Net Load Uncertainty
\thanks{Mr. Yurdakul gratefully acknowledges the support of the German Federal Ministry of Education and Research and the Software Campus program under Grant 01IS17052. }
}

\author{\IEEEauthorblockN{Ogun Yurdakul\IEEEauthorrefmark{1},
Feng Qiu\IEEEauthorrefmark{2}, and
Sahin Albayrak\IEEEauthorrefmark{1}}
\IEEEauthorblockA{\IEEEauthorrefmark{1}Department of Electrical Engineering and Computer Science, Technical University of Berlin, Berlin, Germany}
\IEEEauthorblockA{\IEEEauthorrefmark{2}Energy Systems Division, Argonne National Laboratory, Lemont, IL 60439, USA}}

\maketitle
\begin{abstract}
To take unit commitment (UC) decisions under uncertain net load, most studies utilize a stochastic UC (SUC) model that adopts a \textit{one-size-fits-all} representation of uncertainty. Disregarding contextual information such as weather forecasts and temporal information, these models are typically plagued by a poor out-of-sample performance. To effectively exploit contextual information, in this paper, we formulate a \textit{conditional} SUC problem that is solved \textit{given} a covariate observation. The presented problem relies on the true conditional distribution of net load and so cannot be solved in practice. To approximate its solution, we put forward a predictive prescription framework, which leverages a machine learning model to derive weights that are used in solving a reweighted sample average approximation problem. In contrast with existing predictive prescription frameworks, we manipulate the weights that the learning model delivers based on the specific dataset, present a method to select pertinent covariates, and tune the hyperparameters of the framework based on the out-of-sample cost of its policies. We conduct  extensive numerical studies, which lay out the relative merits of the framework vis-à-vis various benchmarks.   
\end{abstract}
\section{Introduction}\label{1}
Taking unit commitment (UC) decisions under uncertain net load (i.e., load minus renewable generation) lies at the cornerstone of ensuring the economical and reliable operation of systems with deep penetration of renewables. To this end, grid operators (GOs) typically draw upon contextual information (e.g., historical realizations of net load, weather forecasts) as features to train machine learning (ML) algorithms that generate point predictions for net load, which are subsequently utilized in solving a deterministic UC problem. Despite capitalizing on contextual information, such an approach fails to capture the stochastic nature of net load and suffers due to isolating the ML algorithm from the downstream optimization problem. On the flip side, the well-touted stochastic optimization (SO) models explicitly represent uncertainty by usually making assumptions on the probability distribution of net load. Nevertheless, SO models exhibit a poor out-of-sample performance if the assumed probability distribution is wrong, and they cannot effectively exploit covariate observations, resorting to a \textit{one-size-fits-all} representation of uncertainty.  \par
Recently, a paradigm termed predictive prescriptions emerged in the operations research literature, which aims to address these shortcomings by jointly leveraging supervised ML algorithms and a conditional SO model. The approach put forward in \cite{bert} trains an ML model to derive weights for historical observations of the uncertain parameter and uses the weights in solving a reweighted sample average approximation (SAA) problem. Predictive prescription frameworks found applications in power systems as well \cite{pt_ifo, fd_pof, morales}. The study in \cite{pt_ifo} seeks to maximize the profit of a renewable resource trading in the day-ahead market by training trees with a task-based loss, whereas \cite{fd_pof} leverages linear regression models to determine the renewable generation forecasts that lead to UC decisions with minimal total cost. \par
In this paper, we initially map out in Section \ref{sec2} a conditional SO formulation for taking UC decisions under uncertain net load \textit{given} a covariate observation. In Section \ref{sec3}, we put forward a predictive prescription framework, which leverages the random forest (RF) algorithm to derive weights that are used in solving a reweighted SAA problem. Section \ref{sec3} further lays out three principal contributions of this paper. First, we put forth a method that manipulates the weights derived from the RF algorithm based on the size and the information-richness of the dataset. Second, we present an approach to tuning the hyperparameters of the framework based on the out-of-sample cost of its prescriptions. Finally, we suggest a method for pinpointing the pertinent covariates of net load. In Section \ref{sec4}, we demonstrate the application of the framework using data harvested from the California Independent System Operator (CAISO) grid and investigate its out-of-sample and computational performance. Section \ref{sec5} concludes the paper.     
\vspace{-0.0cm}              
\section{Problem Description}\label{sec2}
We start out with the analytical description of the problem. 
\vspace{-0.5cm}
\subsection{Analytical underpinnings}\label{sec2a}
We study the UC problem under the uncertainty in net load, solved by the GO at an hourly granularity for a scheduling horizon of 24 hours. The study period for each day $d$ is denoted by the set $\mathscr{H}_d \coloneqq \{h \colon h=1,...,24\}$, where the term $h$ is the index for each hourly period. We denote by ${Y}_{d} \in \mathscr{Y} \subseteq \mathbb{R}^{d_y}$ the uncertain net load across all system buses and all 24 hours in $\mathscr{H}_d$, and we represent its observation by $Y_d = y_d$. Assume that the GO has at its disposal historical observations on net load for $D$ days. Define the set $\mathscr{D}\coloneqq\{d \colon d=1,\ldots,D\}$.  \par
Typically, it is not possible to precisely set forth the probability distribution of net load ${Y}_{d}$ or provide a perfectly accurate forecast for the materialized net load levels ${Y}_{d} = y_d$. Nevertheless, there is a broad array of contextual information that could prove useful to these ends. For instance, temperature, solar irradiance, and wind speed measurements, temporal information such as the month of the year and the day of the week, as well as lagged observations on net load may have a direct bearing on the net load realization ${Y}_{d} = y_d$. Our framework capitalizes on the observations of these covariates in assessing the uncertainty in net load. Note that such covariates are precisely the features that are leveraged in developing ML models so as to forecast net load. We express by $X_d \in \mathscr{X} \subseteq \mathbb{R}^{d_{x}}$ the random covariate associated with ${Y}_{d}$ and denote its observation by $X_d = x_d$. We expound upon the candidate covariates drawn upon in our framework in Section \ref{sec4}.
\vspace{-0.2cm}
\subsection{Conditional stochastic unit commitment problem}
The contextual information associated with net load can be effectively exploited in a conditional stochastic programming framework in taking UC decisions after observing the contextual information $X= \bar{x}$. If it were possible to know the true, underlying conditional distribution of the net load ${Y}$ given $X = \bar{x}$, we could formulate the following ``gold standard'' conditional stochastic unit commitment ($\mathsf{CSUC}$) problem:
\begin{IEEEeqnarray}{lll}
\hspace{-0.75cm} \underset{z \in \mathscr{Z}}{\text{min}} &\hspace{-0.5cm} \mathbb{E}\big[\mathcal{C}({z};{Y})|{X=\bar{x}} ]\coloneqq&\hspace{.2cm}\eta^{\mathsf{T}}{z} + \mathbb{E}\big[\mathcal{Q}({z};{Y})|{X=\bar{x}}\big]\label{objfs}\\
\hspace{-0.75cm}\text{where} &&\nonumber\\
\hspace{-0.75cm}\mathcal{Q}(z; \bar{y}) \coloneqq & \hspace{0.1cm} \underset{\zeta}{\text{min}}&  \hspace{-0.7cm} c^{\mathsf{T}}\zeta  \label{objss}\\
\hspace{-0.75cm}&  \hspace{0.1cm}\text{subject to} & \hspace{-0.7cm} W\zeta  \leq b - Tz - M\bar{y}.\label{css}
\end{IEEEeqnarray} 
The $\mathsf{CSUC}$ problem is a two-stage conditional SO problem with a mixed-integer linear programming formulation. The objective \eqref{objfs} of the first-stage problem is to minimize the commitment and startup costs plus the expected dispatch and load curtailment costs. The first-stage decisions comprise the binary commitment, startup, and shutdown variables and are represented by $z$. The set $\mathscr{Z}$ denotes the feasible region of the first-stage decisions, which is defined by the logical constraints that relate the commitment, startup, and shutdown variables as well as the minimum uptime and downtime constraints.\par
For a specific vector of first-stage decisions $z$ and materialized net load values $Y = \bar{y}$, the value function $\mathcal{Q}({z};\bar{y})$ is evaluated by solving the second-stage problem \eqref{objss}--\eqref{css} with the objective \eqref{objss} to minimize the dispatch costs and the penalty cost due to load curtailment. The second-stage variables are denoted by $\zeta$ and are composed of the power dispatch levels of generators and the curtailed load for all hours $h \in \mathscr{H}_d$. We succinctly represent in \eqref{css} the power generation and ramping limits for generators as well as the transmission constraints using injection shift factors based on the DC power flow model, wherein $W$, $T$, and $M$ are constant matrices and $b$ is the right-hand-side vector of all second-stage constraints.\par
Note that the $\mathsf{CSUC}$ formulation draws upon the true conditional distribution of $Y | X=\bar{x}$, which cannot be known, thus rendering $\mathsf{CSUC}$ solely a hypothetical, ideal formulation. Nevertheless, GOs have data on the historical realizations of random net load and the associated covariates. The predictive prescription framework introduced in Section \ref{sec3} utilizes these observations to construct the training set  $\mathscr{S}_{D}\coloneqq\{({x}_d, {y}_d)\}_{d=1}^{D}$, which is used to solve a surrogate problem for $\mathsf{CSUC}$.  
\vspace{-0.1cm}
\section{Proposed Framework}\label{sec3}
We next introduce our predictive prescription framework. 
\vspace{-0.1cm}
\subsection{Surrogate problem formulation}\label{sec3a}
The principal objective of the proposed framework is to approximate as close as possible the optimal $\mathsf{CSUC}$ solution after observing $X = \bar{x}$. Motivated by \cite{bert}, we approximate the $\mathsf{CSUC}$ problem by the reweighted SAA problem
\begin{IEEEeqnarray}{l}
\mathsf{w-CSUC}\colon\,\hat{z}_{D}(\bar{x})\in\text{arg}\,\underset{z \in \mathscr{Z}}{\text{min}}\hspace{0.2cm}\sum_{d=1}^{D}{\omega}_{D, d}(\bar{x}) \mathcal{C}({z};{y_d}),\label{objwSUC}
\vspace{-.1cm}
\end{IEEEeqnarray}
where ${\omega}_{D, d}(\bar{x})$, $\forall d \in \mathscr{D}$, are weight functions obtained from the training set that adjust the influence of each historical observation on the objective function of the reweighted SAA problem \eqref{objwSUC}. In \cite{bert}, the authors derive the weights ${\omega}_{D, d}(\bar{x})$ by directly using ML algorithms and spell out how different ML algorithms can be leveraged to that end. 
In the proposed framework, we draw upon the RF model in conjunction with a nonlinear function to derive the weights. 
\vspace{-0.1cm}
\subsection{Evaluation of the empirical weights}\label{sec3b}
We start off by training the RF model to predict the net load in the next 24 hours, for which we use the covariate observations in $\mathscr{S}_{D}$ as features and the net load values in $\mathscr{S}_{D}$ as labels. Next, for each new covariate observation $X=\bar{x}$, we use the trained RF model to quantify the similarity between the new observation and each historical observation in $\mathscr{S}_{D}$.\par
To quantify the similarity between observations, we record the leaf that the new observation $\bar{x}$ is mapped into in each tree of the RF and subsequently identify the historical covariate observations that fall into the same leaf node with $\bar{x}$. Central to our approach is to assign the weight for observation $x_d$, that is, ${\omega}_{D, d}(\bar{x})$, based on the number of trees in which $x_d$ and $\bar{x}$ are assigned to the same leaf node. To this end, Bertsimas and Kallus \cite{bert} propose that the weight for $x_d$ increase linearly with the number of trees in which $\bar{x}$ and $x_d$ fall into the same leaf node, normalized by the total number of covariate observations assigned to the same leaf node with $\bar{x}$, which yields the \textit{empirical} weights
\begin{IEEEeqnarray}{c}
\hat{\omega}_{D, d}(\bar{x}) =\frac{1}{T}\sum_{\tau=1}^{T}\frac{\mathbb{I} \big[x_d \in \mathcal{X}^{\tau}_{l(\bar{x})}\big]}{\big|\big\{d' \colon x_{d'} \in \mathcal{X}^{\tau}_{l(\bar{x})} \big\}\big|},
\end{IEEEeqnarray}
where $T$ denotes the number of trees in the forest, $\mathbb{I}(\cdot)$ the indicator function, and $\mathcal{X}^{\tau}_{l(\bar{x})}$ the set of covariate observations assigned to the same leaf with $\bar{x}$.\par
\subsection{Deriving the final weights}\label{sec3c}
Existing approaches in the literature plug the empirical weights $\hat{\omega}_{D, d}(\bar{x})$ into the $\mathsf{w-CSUC}$ problem without taking into account the size or the prescriptive content of the training set. Nevertheless, the empirical weights obtained with a small training set and/or a training set with little informative content may fail to afford an accurate characterization of the similarity between observations. In contrast, we can utilize a particular historical observation with greater confidence, if it were deemed to be similar to a new observation under a large training set with high prescriptive power. 
As such, we introduce the function $\varphi\big(\hat{\omega}_{D, d}(\bar{x}); \xi, D \big)\coloneqq$ $\hat{\omega}_{D, d}(\bar{x}) ^{\frac{D}{\xi}}$, which serves to manipulate the empirical weights based on the training set size $D$ and the weight modification parameter $\xi$. The parameter $\xi$ could be viewed as a proxy for the information richness and the prescriptive power of the training set. As $\xi$ decreases and $D$ increases, $\varphi(\cdot)$ amplifies the weights of the data points that are assessed to be strongly similar to a new observation and brings down the empirical weights of the points that are markedly dissimilar to a new observation. Further, for a small $D$ and a large $\xi$, $\varphi(\cdot)$ smoothens any significantly high and low weight value and brings the weights toward a uniform level. 
Clearly, a key challenge to this end is to hone in on a judicious value of $\xi$.
To this end, we treat $\xi$ as a hyperparameter of the overall framework and set its value by assessing its influence on a separate validation set. 
\vspace{-0.25cm}
\subsection{Task-based hyperparameter tuning}\label{sec3d}
The proper tuning of an ML model's hyperparameters could play a drastic role in its performance. The classical approach to hyperparameter tuning is to assess the performance of an ML under different hyperparameter values based on a statistical loss function. In our framework, however, a specific selection of RF hyperparameter values may bring forth a lower prediction error without leading to UC decisions that drive down the total out-of-sample cost, which punctuates the need to tune the RF model's hyperparameters based on the ultimate task for which it is trained. As such, we treat the hyperparameters of the RF model as those of the overall framework and set their values based on the total out-of-sample cost of the optimal policy $\hat{z}_{D}(\bar{x})$ obtained with different hyperparameter values. \par
At the outset, we use grid search to exhaustively generate candidate values for the hyperparameters reported in Table \ref{hyperpa}. Next, we construct a separate validation set containing pairs of covariate and net load observations $\mathscr{V}_{\bar{D}}\coloneqq\{(\bar{x}_i,\bar{y}_i)\}_{i=1}^{\bar{D}}$. We use each covariate observation $\bar{x}_i$ to compute the optimal policy $\hat{z}_{D}(\bar{x}_i)$ and subsequently the out-of-sample cost $\mathcal{C}(\hat{z}_{D}(\bar{x}_{i}); \bar{y}_{i})$ obtained under the actual net load observation $\bar{y}_i$. For each set of candidate hyperparameter values, we compute the total out-of-sample cost over the validation set, i.e., $\sum_{i=1}^{\bar{D}} \mathcal{C}(\hat{z}_{D}(\bar{x}_{i}); \bar{y}_{i})$. Ultimately, we pick the hyperparameter values that deliver the lowest total out-of-sample cost. 
\begin{table}[h]
\vspace{-.25cm}
\centering
{\fontsize{8.5}{11.05}\selectfont
\setlength{\tabcolsep}{7pt} 
\renewcommand{\arraystretch}{1.2}
\caption{Hyperparameters}
\label{hyperpa}
\centering
\begin{tabular}{c | c  }
 \hline \hline 
{hyperparameter} & candidate values\\
\hline \hline
max tree depth & 3, 6, 10\\ \hline
number of features considered  & \multirow{2}{*}{$\sqrt{d_{x}}$, $(0.3){d_{x}}$, $(0.6){d_{x}}$} \\ 
for node splitting & \\ \hline
weight modification parameter $\xi$ & $\frac{D}{10}$, $\frac{D}{4}$, ${D}$, ${4D}$, $10D$\\
\hline \hline
\end{tabular}}
\vspace{-.5cm}
\end{table}
\vspace{0cm}
\subsection{Selection of the covariates}\label{sec3e}
A key thrust of the framework is to pinpoint information-rich covariates that aid in effectively grasping the uncertainty in net load. While there is plethora of factors that can potentially influence a net load realization, ruling out the covariates that afford little or no information can help the trained RF model better assess the similarity between covariate observations, thereby yielding final weights that reduce the out-of-sample costs. Further, working with fewer covariates allows for expediting the training and testing of the RF model.\par
To identify the covariates, we employ a hybrid approach comprising a filter and a wrapper feature selection method. We denote the support of the initial set of candidate covariates by $\mathscr{X}^{r} \subseteq \mathbb{R}^{d_{x^{r}}}$. We start out by computing the Pearson correlation coefficient (PCC) for each candidate covariate and the net load observation, which measures the linear correlation between two variables. The PCC attains values between $-1$ and $1$, with $1$ (resp.  $-1$) indicating a complete positive (resp. negative) correlation and $0$ signifying that the correlation is immaterial. Customarily, when the absolute value of the PCC is greater than or equal to $0.6$, it is interpreted as the variables being strongly correlated with one another \cite{pcc_num}. As such, we rule out all candidate features that yield a PCC value between $0.6$ and $-0.6$ and obtain $d_{x^{p}}$ covariates supported on the set $\mathscr{X}^{p} \subseteq \mathbb{R}^{d_{x^{p}}}$. Note that PCC measures only the linear correlation between the variables, and it does not assess how the covariates integrate with the utilized ML model. As a remedy, we additionally implement recursive feature elimination (RFE), which is a wrapper method that takes an ML model as a parameter. RFE trains the selected ML model with the initial set of features, ranks the features on the basis of their importance, and recursively eliminates the least important features until the desired number of features is reached. We run RFE with the RF model 
and ultimately obtain the final set of covariates with support $\mathscr{X}^{f} \subseteq \mathbb{R}^{d_{x^{f}}}$. 
\section{Numerical Experiments}\label{sec4}
We next demonstrate the application of the proposed framework in a real-life setting.
\subsection{Datasets and covariate selection}\label{sec4a}
In our experiments, we draw upon the net load values recorded in the CAISO grid between June 1, 2018 and August 31, 2019 \cite{ols:caiso}. We use the measurements recorded in the first year (i.e., June 1, 2018--May 31, 2019) to construct the training sets. As set forth below, we vary the size of the training set in different experiments so as to assess its influence on the performance of the methods. Nevertheless, to compare their performance on a consistent basis, we use the same validation set and the same test set in all experiments. Specifically, we utilize the measurements recorded in June 2019 as the validation set and the measurements recorded from July 1 to August 31, 2019 as the test set. We use the IEEE 14-bus system 
in the experiments, which has 5 generators with an aggregate capacity of 765.31 MW. We scale the net load values so that the highest net load value is equal to the 90\% of the aggregate capacity of the generators.\par
To select the covariates for PV and wind generation, we assess the spatial distribution of PV and wind installations with their respective capacities across California, and we accordingly select locations from which we harvest data on global horizontal irradiance (GHI) and wind speed (magnitude and direction). 
We study the total population and population density of the counties of California and identify locations from which we use temperature measurements so as to capture the influence of temperature on system load.\footnote{We provide the data and the source code of the simulations in the online companion to this paper located in \url{https://github.com/oyurdakul/isgtna23}.} We leverage as candidate covariates the GHI, wind speed, and temperature measurements reported by the National Renewable Energy Laboratory \cite{nrel} for the selected locations in the past 24 hours. We further use as candidate covariates 24 lagged realizations of net load, as well as the 24 lagged realizations of the daily, weekly, and monthly moving average of net load. Finally, we define categorical variables to indicate whether a day falls on a weekend and on a public holiday and use one-hot encoding for their representation. We ultimately obtain $d_ {x^{r}}= 440$ candidate covariates and follow the covariate selection method presented in Section \ref{sec3d} to derive $d_{x^{f}}=25$ covariates.
\subsection{Benchmarks}\label{sec4b}
To highlight the relative merits of the proposed framework, we draw upon different decision-making methods to obtain alternative policies and investigate their performance. One such method is the naive stochastic unit commitment ($\mathsf{NSUC}$) model, which treats the net load observations in the training set as equiprobable scenarios and disregards the covariate observation $X = \bar{x}$, stated as
\begin{IEEEeqnarray}{lll}
\mathsf{NSUC:}\hspace{0.25cm} \underset{z \in \mathscr{Z}}{\text{min}} &\hspace{0.25cm} \frac{1}{D}\sum_{d=1}^{D}\mathcal{C}({z};{y_{d}}). \label{objnsuc}
\end{IEEEeqnarray}
We solve the following reweighted SAA problem using the empirical weights as suggested in \cite{bert} so to investigate the impact of the transforming the weights:
\begin{IEEEeqnarray}{lll}
\mathsf{ew-CSUC:}\hspace{0.25cm} \underset{z \in \mathscr{Z}}{\text{min}} &\hspace{0.25cm} \sum_{d=1}^{D}{\hat{\omega}}_{D, d}(\bar{x}) \mathcal{C}({z};{y_d}). \label{objewuc}
\end{IEEEeqnarray}
We further use the point forecast of the trained RF model, i.e., $\hat{f}_{D}^{RF}(\bar{x})$, in solving the following deterministic UC problem:
\begin{IEEEeqnarray}{lll}
\mathsf{PFUC:}\hspace{0.25cm} \underset{z \in \mathscr{Z}}{\text{min}} &\hspace{0.25cm} \mathcal{C}({z};\hat{f}_{D}^{RF}(\bar{x})). \label{objpfuc}
\end{IEEEeqnarray}
To obtain the minimum out-of-sample cost that could be ideally attained, we solve the following ideal UC ($\mathsf{IUC}$) problem, which has a perfect foresight of the net load observation $\bar{y}$:   
\begin{IEEEeqnarray}{lll}
\mathsf{IUC:}\hspace{0.25cm} \underset{z \in \mathscr{Z}}{\text{min}} &\hspace{0.25cm} \mathcal{C}({z};\bar{y}). \label{objiuc}
\end{IEEEeqnarray}
\vspace{-0.4cm}
\subsection{Results}
We conduct the experiments on a 64 GB-RAM computer containing an Apple M1 Max chip with 10-core CPU. We build the RF models and select the covariates under Python using scikit-learn 1.1.2. To model the UC instances, we extend the \verb|UnitCommitment.jl| package \cite{UCjl} to the two-stage stochastic setting, and we solve all UC problems under Julia 1.6.1 with Gurobi 9.5.0 as the solver. The penalty cost for load curtailment is set at $\$10,000/MWh$ in all experiments. \par
We initially construct the training set with the first 100 observations, i.e., $D=100$. We use the validation set to tune the hyperparameters of the proposed predictive prescription framework (indicated as $\mathsf{w-CSUC}$) as well as those of the $\mathsf{PFUC}$ and $\mathsf{ew-CSUC}$ methods. To assess how the methods perform out-of-sample, we use the measurements in the test set to determine how each method would have committed the generators for the corresponding days in the test set, then we observe the actual net load levels that had materialized, and ultimately use the resulting total cost and the mean unserved energy (MUE) to score the performance of each method. In Table \ref{oos_100}, we report for each method the average of the total cost and the MUE computed over all 62 observations in the test set. We separately tabulate the MUE results in addition to the total cost as the latter may greatly vary with the choice of the penalty cost for load curtailment. \par
\begin{table}[h]
\vspace{-.3cm}
\centering
{\fontsize{8.5}{11.05}\selectfont
\setlength{\tabcolsep}{7pt} 
\renewcommand{\arraystretch}{1.2}
\caption{Out-of-sample costs and MUE levels}
\label{oos_100}
\centering
\begin{tabular}{c | c  | c}
 \hline \hline 
{method} & total cost (\$) & MUE $(MWh)$ \\
\hline \hline
$\mathsf{IUC}$ & 380089.8 & 0.0\\ \hline
$\mathsf{w-CSUC}$ & 401377.3 & 0.0\\ \hline
$\mathsf{ew-CSUC}$ & 414566.0 & 0.0\\ \hline
$\mathsf{NSUC}$ & 416429.5 & 0.0\\ \hline
$\mathsf{PFUC}$ & 453638.1 & 7.2\\
\hline \hline
\end{tabular}}
\vspace{-.3cm}
\end{table}
\begin{figure*}%
    \hspace{-0.24cm}
    \subfloat[][\label{res_1a}]{{\includegraphics[width=.46\linewidth]{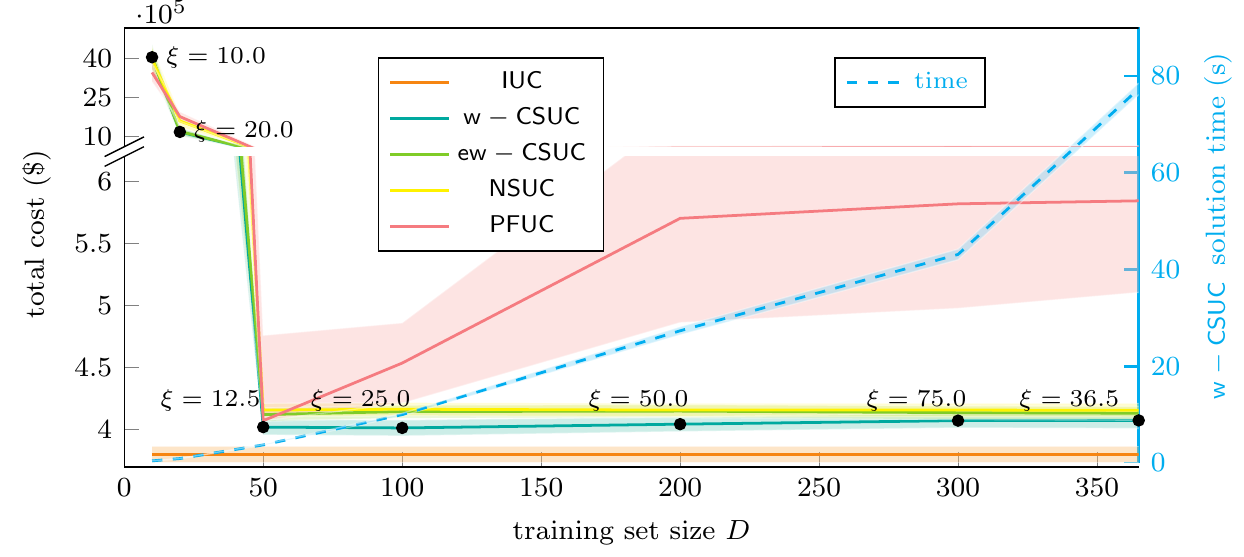} }}\hspace{1.27cm}%
    \subfloat[][\label{res_1b}]{{\includegraphics[width=.46\linewidth]{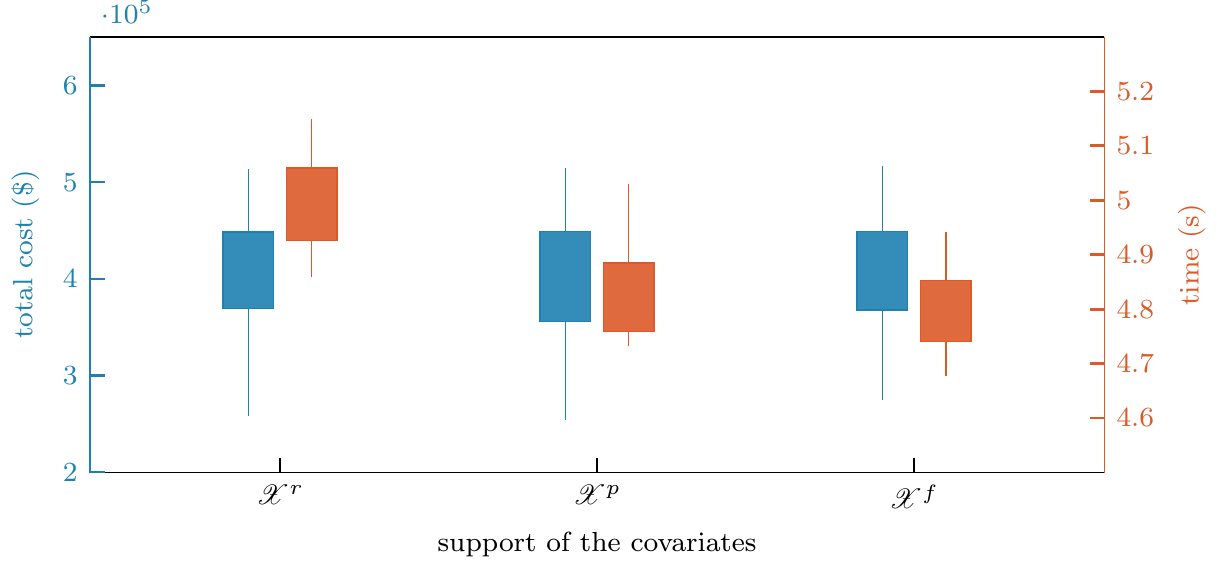} }}%
    \vspace{-.23cm}
\caption{Out-of-sample performances. In \ref{res_1a}, solid lines indicate the mean values over the 62 observations in the test set, and shaded regions show only one-tenth of the standard deviations around the mean values in order to avoid excessive overlap between shaded regions. \ref{res_1b} depicts the total cost and the computation time for obtaining  the weights ${\omega}_{D, d}(\bar{x})$ under each set of features. }
    \label{test}%
    \vspace{-0.55cm}
\end{figure*}
The results in Table \ref{oos_100} make clear the monetary benefits that can be reaped by implementing $\mathsf{w-CSUC}$, which delivers the lowest total cost among all methods except the perfect-foresight policy, yielding a total cost that is higher by 3.86\% than IUC. We highlight that the lower cost under $\mathsf{w-CSUC}$ in comparison with that under the runner-up method $\mathsf{ew-CSUC}$ provides an empirical justification for manipulating the empirical weights before using them in solving the reweighted SAA problem. The $\mathsf{NSUC}$ method fails to outperform $\mathsf{w-CSUC}$ and comes at the heels of $\mathsf{ew-CSUC}$, which we ascribe to $\mathsf{NSUC}$ utilizing equiprobable scenarios without taking the covariate observations into account.  The results further make evident the shortcomings of drawing upon deterministic forecasts and ignoring the stochastic nature of net load in solving the UC problem, as the policies under $\mathsf{PFUC}$ deliver the highest total cost and necessitate involuntary load curtailment.\par
In certain practical applications, we may fail to collect a large number of observations that can be used in constructing the training set. Coincidentally, working with a larger training set requires solving the $\mathsf{w-CSUC}$ problem with a greater number of scenarios, which drives up the computational burden. As such, we assess the performance of each method under different values of $D$. In doing so, we keep all the hyperparameters except $\xi$ constant at the values determined for $D=100$ and use grid search to tune the value of $\xi$ on the validation set. Fig. \ref{res_1a} visualizes the total cost delivered by each method, and it illustrates for the proposed framework the value of $\xi$ determined through grid search and the time to solve the $\mathsf{w-CSUC}$ problem so as to obtain the policy $\hat{z}_{D}(\bar{x})$. \par
The plots in Fig. \ref{res_1a} echo the order of performance in Table \ref{oos_100}, as across most values of $D$, the policies of $\mathsf{w-CSUC}$ beat $\mathsf{ew-CSUC}$, which in turn outperforms $\mathsf{NSUC}$. Note that the policies under $\mathsf{PFUC}$ exhibit the worst out-of-sample performance for most investigated training set sizes. We further point out that the total costs obtained under the proposed framework are tightly clustered around their mean values and less spread out compared with the benchmark methods.\par
We remark upon the tight coupling between the training set size, the solution time, and the out-of-sample cost. Increasing $D$ from $10$ to $100$ markedly improves (decreases) the out-of-sample performance of $\mathsf{w-CSUC}$, albeit a diminishing return as $D$ grows from $50$ to $100$. We also observe that the out-of-sample performance saturates around $D=100$ and sporadically deteriorates (increases) as $D$ grows beyond $100$, during which the solution time precipitously increases. \par
One can draw from Fig. \ref{res_1a} valuable insights into the value of $\xi$ determined via grid search. Most notably, the empirical weights are amplified and suppressed the most under $D=365$, signifying an information-rich training set. This observation drives home that training the RF model using a full year’s data enables an accurate characterization of the similarity between observations. Note that, the empirical weights for $D \in \{50, 100, 200, 300\}$ are also boosted and attenuated, though not as much as for $D=365$, whereas those for $D \in \{10, 20\}$ are used as is, indicating that amplifying and suppressing the empirical weights obtained with such small datasets is not warranted. \par
We next investigate the influence of the covariate selection method laid out in Section \ref{sec3e} on the out-of-sample performance and the computation time. To this end, we repeat the experiments for $D=100$ under the set covariates supported in $\mathscr{X}^r$ and $\mathscr{X}^p$. We tune the hyperparameters for each set of covariates using the validation set and leverage the test to compute the out-of-sample performances. For each set of covariates, we measure the time for computing the weights ${\omega}_{D, d}(\bar{x})$ over 30 simulation runs, which is comprised of the time for training the RF model and that for evaluating the weight for all observations in the test set. Fig. \ref{res_1b} bears out the relative merits of the proposed covariate selection method, which notches a 3.90\% reduction in the average time for evaluating the weights ${\omega}_{D, d}(\bar{x})$ vis-à-vis those under the initial set of features without compromising on the out-of-sample performance.
\vspace{-.1cm} 
\section{Conclusion}\label{sec5}
In this paper, we worked out a predictive prescription framework that jointly leverages the random forest (RF) algorithm with a conditional stochastic optimization model so as to take unit commitment decisions under uncertain net load. We put forth a method to manipulate the empirical weights derived from the RF model based on the size and the prescriptive power of the training set, and we suggest a hybrid method to select pertinent covariates for net load. By treating the hyperparameters of the RF model as those of the overall framework, we tune them based on the ultimate task for which the framework is developed, that is, bringing forth a lower out-of-sample cost. The extensive numerical studies conducted illustrate the capabilities of the framework in reducing not only the out-of-sample cost and load curtailment, but also the computation time compared with various benchmarks.  
\vspace{-.1cm}
\bibliographystyle{IEEEtran}
\bibliography{IEEEabrv,csobib}
\end{document}